%



\vsize=230mm
\font\titre=cmbx10 scaled \magstep1

\def\carre{\kern 0.12pt\vcenter{\hrule \hbox{\vrule height 3pt
           \kern 3pt \vrule}\hrule}\kern 0.24pt}
\def\Carre{\kern 0.16pt\vcenter{\hrule \hbox{\vrule height 4pt
           \kern 4pt \vrule}\hrule}\kern 0.32pt}

\def\cfBo{[\kern 1.5pt  1\kern 1.5pt]}
\def\cfCh{[\kern 1.5pt  2\kern 1.5pt]}
\def\cfEh{[\kern 1.5pt  3\kern 1.5pt]}
\def\cfFLM{[\kern 1.5pt 4\kern 1.5pt]}
\def\cfJS{[\kern 1.5pt  5\kern 1.5pt]}
\def\cfLe{[\kern 1.5pt  6\kern 1.5pt]}
\def\cfMi{[\kern 1.5pt  7\kern 1.5pt]}
\def\cfMS{[\kern 1.5pt  8\kern 1.5pt]}
\def\cfPi{[\kern 1.5pt  9\kern 1.5pt]}
\def\cfPib{[\kern 1.5pt 10\kern 1.5pt]}
\def\cfPr{[\kern 1.5pt  11\kern 1.5pt]}
\def\cfTa{[\kern 1.5pt  12\kern 1.5pt]}
\def\cfTab{[\kern 1.5pt 13\kern 1.5pt]}

\centerline{\titre Some deviation inequalities}
\bigskip
\centerline{by Bernard Maurey}
\medskip
\centerline{\sevenrm September 1990}
\bigskip

\bigskip

\noindent {\sevenrm {\sevenbf Abstract.} We introduce a
concentration property for
probability measures on $\scriptstyle{R^n}$, which we call
Property~($\scriptstyle\tau$);
we show that this property has an interesting stability under
products and contractions
(Lemmas 1,~2,~3). Using property~($\scriptstyle\tau$),
we give a short proof for a recent
deviation inequality due to Talagrand. In a third
section, we also recover
known concentration results for Gaussian measures using our approach.}

\bigskip

\noindent{\bf Introduction.} Very roughly speaking, a concentration
of measure phenomenon for a probability measure~$\mu$ means that,
given any
measurable subset~$A$ such that $\mu(A) \ge {1\over 2}$, the
enlargements of~$A$, in a sense to be made precise, almost have
measure~1. One important example of this situation is a
consequence of Paul L\'evy's isoperimetric inequality for the unit
sphere $S_{n}$ of~$R^{n+1}$. This consequence is the following: Let
$\mu_n$ denote the normalized rotation invariant measure on~$S_{n}$;
for every measurable subset $A$ of~$S_{n}$ such that
$\mu_n(A) \ge {1\over 2}$ we have
$$ \mu_n\{ x\in S_n; x\notin A_\varepsilon \} \le
\sqrt{\pi\over 8} \, e^{-n\varepsilon^2/2}$$
where $A_\varepsilon$ denotes the subset of $S_{n}$ of all points
whose geodesic distance to~$A$ is less than~$\varepsilon$
(see~\cfMS).
This result is crucial for the proof of Dvoretzky's theorem
(about almost spherical sections of convex bodies) as given
by V.~Milman~\cfMi\
(see also Figiel, Lindenstrauss and Milman~\cfFLM).
Using L\'evy's result, Borell~\cfBo\ was able to prove an analogous
isoperimetric result for the Gaussian measure on~$R^n$
(see also Ehrhard~\cfEh).
As it is the case for the sphere, Borell's isoperimetric result
implies a Gaussian concentration of measure principle.
Later, Pisier and the author gave a very
simple proof for the Gaussian concentration principle which is needed
for the proof of Dvoretzky's theorem (see~\cfPi,~\cfMS\ and~\cfPib).
Recently, Talagrand proved~\cfTa\ a concentration principle
for measures on~$R^n$
with exponential densities
which is stronger than the Gaussian
one (Corollary~1 below; this Corollary appears in~\cfTa\ as a
consequence of a more precise isoperimetric inequality
which does not follow from our proof).
  We give here a proof for Talagrand's concentration result, using a
property which we call Property~($\tau$). This property is defined in
section~I, and general stability results
about it are stated (Lemmas~1, 2,~3); we also explain in section~I
why property~($\tau$) is related to concentration (Lemma~4).
In section~II, we show
how to get Talagrand's inequality using property~($\tau$). In
section~III we recover the Gaussian case along the same
lines (Corollary~2).
In a last section, we study a variant of property~($\tau$),
the {\it convex} property ($\tau$); this variant is related to
an other deviation inequality due to Talagrand~\cfTab\ (see
Corollary~5).
\bigbreak

\noindent{\bf I. Property ($\tau$)}
\medskip

 Let $f$ and $g$ be two measurable functions on~$R^n$; we denote
by $f \Carre g$ the inf-convolution of $f$ and~$g$,
$$ (f \Carre g)(x) = \inf \{ f(x-y) + g(y);y\in R^n \}.$$
 If $\mu$ is a probability measure on~$R^n$ and~$w$ a positive
measurable function on~$R^n$, we say that the couple $(\mu, w)$
satisfies the Property~($\tau$) if for every bounded measurable
function~$\varphi$ on~$R^n$ we have
$$ (\int e^{ \varphi\carre w} d\mu)(\int e^{ -\varphi} d\mu)
\leq 1.$$
   If we adopt the convention $+\infty.0 \leq 1$, we see easily that the
above inequality extends to all $\overline R$-valued
measurable functions $\varphi$.

  The definition of property~($\tau$) was motivated by
Talagrand's isoperimetric inequality for the cube, as
the careful reader of \cfTab\ will notice.
\smallbreak

\proclaim Lemma 1. If $(\mu_i, w_i)$ satisfies~($\tau$) on~$R^{n_i}$
for~$i=1,2$,
then $(\mu_1\otimes \mu_2, w)$ satisfies~($\tau$) on~$R^{n_1}\times
R^{n_2}$, where
$$ w(x_1,x_2) = w_1(x_1) + w_2(x_2).$$

\noindent Proof: Consider $\varphi^y(x)=\varphi(x,y)$ and apply
($\tau$) to
$\psi(y)=\log(\int e^{\varphi^y\carre w_1} d\mu_1).$
\medbreak

\proclaim Lemma 2. If $(\mu_i, w_i)$ satisfies~($\tau$) on~$R^{n}$
for~$i=1,2$,
then $(\mu_1 * \mu_2, w_1 \Carre w_2)$ satisfies~($\tau$) on~$R^n$.
\medbreak

\proclaim Lemma 3. Let $(\mu_1, w_1)$ satisfy~($\tau$) on~$R^{n_1}$. Let
$w_2$ be a positive measurable function
on $R^{n_2}$ and $F$ a mapping from $R^{n_1}$ to~$R^{n_2}$ such that
$w_2(Fx-Fy) \leq w_1(x-y)$ for every pair $x,y$.
Let $\mu_2$ be the image
probability measure on~$R^{n_2}$ defined by $\mu_2=F(\mu_1)$.
Then $(\mu_2,w_2)$ satisfies~($\tau$).

\noindent Proof: $(\varphi \circ F) \Carre w_1 \geq
(\varphi \Carre w_2) \circ F$.
\medbreak
 Lipschitz maps like the above~$F$ were used by
Pisier~\cfPi\ in a
slightly different context. Following him, we will be able in section~III
to pass from the Gaussian case to the uniform probability measure
on $[0,1]^n$ (Corollary~4).

\bigbreak
  We have perhaps to explain why property~($\tau$) for a
couple $(\mu,w)$ is a concentration of measure property.

\proclaim Lemma 4. Assume that $(\mu,w)$ satisfies~($\tau$) on~$R^n$.
For
every measurable subset $A$ of~$R^n$ and every positive real
number~$t$, we have
$$ \mu\{x \notin A + \{w<t\} \} \le (\mu(A))^{-1} e^{-t}.$$

\noindent Proof:
 Let $A$ be a measurable subset of~$R^n$ and denote by $\varphi_A$
the function equal to $0$ on~$A$ and $+\infty$ outside.
We observe that $(\varphi_A\Carre w)(x) \ge t$ when
$x \notin A + \{w<t\} = \{a+y; a\in A, w(y)<t \}$.
Property~($\tau$) implies that
$\int e^{\varphi_A\carre w} d\mu \le (\mu(A))^{-1}$ and we
conclude with Tschebycheff's inequality.
\bigbreak

\noindent{\bf II. Talagrand's deviation inequality}
\medskip

 Let us define a function $W$ on $R$ by
$$ W(t) =  {1\over 18}t^2 \hbox{\rm\quad for\quad} |t|\leq 2,
   \quad{2\over 9}(|t| - 1) \hbox{\rm\quad otherwise}$$
and let $\mu_e$ be the probability measure on $R$ with density
$1_{(0,\infty)}(x) e^{-x}$.

\medskip
\proclaim Proposition. The couple $(\mu_e, W)$ satisfies ($\tau$).
\smallskip

 It follows from Lemma~2 that $(\xi, U)$ also satisfies~($\tau$),
where $\xi$ is the convolution of $\mu_e$ and its symmetric image
on $(-\infty,0)$, and $U = W \Carre W$.
 It is easy to see that $\xi$ has  density ${1\over 2} e^{-|x|}$
on~$R$, and that
$$U(t)=2W(t/2)= {1\over 36} t^2 \hbox{\rm\quad for\quad} |t|\leq 4,
\quad{2\over 9}(|t| - 2)
\hbox{\rm\quad otherwise}.$$
   We deduce now from Lemma~1 that the couple $(\xi_{n}, U_n)$
satisfies~($\tau$) on~$R^n$ for every~$n$, where $\xi_{n}$ is the
product of $n$ copies of $\xi$ and $U_n(x)=\sum_{i=1}^n U(x_i)$.
The idea of working with functions like $W$ or~$U_n$ comes from
Talagrand~\cfTa.
\smallbreak

\proclaim Theorem 1. The couple $(\xi_n, U_n)$ satisfies ($\tau$)
for every integer~$n$. In particular,
for every measurable subset $A$ of~$R^n$ we have
if we set $\rho_A = \varphi_A\Carre U_n$
$$ \int e^{\rho_A} d\xi_{n}  \leq  (\xi_{n}(A))^{-1}.$$

\medbreak
\proclaim Corollary 1 (Talagrand). For every $t>0$,
$$ \xi_n \{x; x\notin A+ 6\sqrt{t} B_2 + 9tB_1 \} \leq
 (\xi_n(A))^{-1} e^{-t}.$$
where $B_2$ and $B_1$ are respectively the usual $\ell_2^n$ and
$\ell_1^n$ balls.

\noindent Proof: According to Lemma~4, we need only show that
$$ \{U_n < t \} \subset 6\sqrt{t} B_2 + 9t B_1.$$
Assume $U_n(x) < t$, and define $y$ and~$z$ in the following way:
$y_i = x_i$ if $|x_i | \leq 4$, $y_i = 0$ otherwise; $z_i = x_i$
if $|x_i | > 4$, $z_i = 0$ otherwise. Then $x = y + z$
and it is easy to check that
$\|y\|_2 \leq 6\sqrt{t}$, $\|z\|_1 \leq 9t$.
\bigbreak

 We present now the proof of the above Proposition.
Let~$\varphi$ be a bounded measurable function on $(0,+\infty)$,
and let~$\psi$ denote the function $\varphi \Carre W$.
\smallskip
Let $I_0 = \int_0^\infty e^{-\varphi(x) -x} dx$ and
$I_1 = \int_0^\infty e^{\psi(y)-y} dy$.
For~$t\in (0,1)$, we define $x(t)$ and~$y(t)$ by the relations
$$ \int_0^{x(t)} e^{-\varphi(x)-x} dx = t I_0,
\int_0^{y(t)} e^{\psi(y)-y} dy = t I_1.$$

We obtain by differentiation
$$ x'(t)=I_0 e^{\varphi(x(t))+x(t)}, y'(t)=I_1 e^{-\psi(y(t))+y(t)}.$$

 Taking into account the
fact that $\psi(y(t)) \leq \varphi(x(t)) + W(x(t)-y(t))$, we obtain
$$ y'(t) \geq I_1 e^{-\varphi(x(t))-W(x(t)-y(t))+y(t)}.$$

Let now $z(t)= {1\over 2} (x(t)+y(t)) - W(x(t)-y(t))$.
We have
 $$z'(t) = \bigl({1\over 2} - W'(x(t)-y(t))\bigr)x'(t) +
 \bigl({1\over 2} + W'(x(t)-y(t))\bigr)y'(t).$$
 If we use the fact that $|W'| \leq 1/2$ on $R$, we get
(writing $x$ and $y$ for $x(t)$ and $y(t)$ for simplicity, and using
the inequality ${1\over 2}(ua+v/a) \geq \sqrt{uv}$
with $a=e^{\varphi(x)}$)

$$ z'(t) \geq (1 - 2W'(x-y))I_0 e^{x} \, {e^{\varphi(x)} \over 2}
+ (1+ 2W'(x-y))I_1 e^{-W(x-y)+y} \, {e^{-\varphi(x)} \over 2}$$
$$ \geq \sqrt{1-4W'(x-y)^2} \, \sqrt{I_0 I_1} \,
e^{{1\over 2}(x+y) - {1\over 2}W(x-y)}$$
$$ = \sqrt{I_0 I_1} \, e^{z(t)} \sqrt{1 - 4W'(x-y)^2} \,
e^{{1\over 2}W(x-y)}.$$
   We claim that for every~$s$
$$ (1 - 4W'(s)^2) \, e^{W(s)} \geq 1.$$
It will then follow that $e^{-z(t)} z'(t) \geq \sqrt{I_0 I_1}$, which
yields after integrating between $0$ and~$1$
$$ 1 \geq \sqrt{I_0 I_1}$$
and this is our Proposition.
\medbreak

\noindent Proof of the claim: We only consider $s \geq 0$
since $W$ is even. For
$s \geq 2$, $W'$~is constant and $W$ increasing, so it is enough to
check the case $0 \leq s \leq 2$; this reduces to
$$ e^{-u/18} \leq 1- 4u/81 \quad\hbox{\rm  for }\quad u \in (0,4)$$
which is proved using elementary calculus.
\bigbreak

\noindent{\bf III. The Gaussian case}
\medskip

   Let $\gamma$ be the standard Gaussian probability measure on~$R$, with
density ${1 \over \sqrt{2 \pi}} e^{-x^2/2}$, and~$\gamma_n$ the product
of $n$~copies of $\gamma$.
   Throughout this section, the norm will be the Euclidean norm
on~$R^n$.

\proclaim Theorem 2. The couple $(\gamma_n, {1\over 4}\|x\|^2)$
satisfies~($\tau$) for every integer~$n$.

\noindent Proof: We check first that $(\gamma,x^2/4)$ has
property~($\tau$)
on~$R$; the proof is similar to the proof of the Proposition,
but simpler: $x(t)$ and~$y(t)$ are defined in a similar fashion,
and $z(t)$ is simply equal to ${1\over 2} (x(t)+y(t))$.
It follows from Lemma~1 that $(\gamma_n, {1\over 4}\|x\|^2)$
has property~($\tau$) for every integer~$n$.

  We can also give a direct proof using the functional
Brunn-Minkowski inequality due to Prekopa and Leindler~\cfPr,
\cfLe\ (see also \cfPib):
If $f,g,h$ are bounded below measurable functions on~$R^n$ such that
for all~$x$ and~$u$ we have
${1\over 2}\bigl(f(x+u)+g(x-u)\bigr) \ge h(x) $, then
$$ (\int e^{-f(x)} dx)(\int e^{-g(x)} dx) \le (\int e^{-h(x)} dx)^2.$$
We apply this inequality to $f(x)= \varphi(x) + {1\over 2}\|x\|^2,
g(y) = -\psi(y) + {1\over 2}\|y\|^2$ and $h(z) = {1\over 2}\|z\|^2$,
where we have set $\psi = \varphi \Carre w,
w(y) = {1\over 4} \|y\|^2$.

\medbreak
\noindent{\bf Remark 1.} The second proof of Theorem~2 only uses
the uniform convexity
properties of $-\log f$, where~$f$ is the density of $\mu$.

\medbreak
\noindent{\bf Remark 2.} As pointed out by Talagrand, the Gaussian
concentration result is a consequence of Corollary~1, using a suitable
map that transforms $\xi_n$ into~$\gamma_n$. More precisely, Lemma~3
and Theorem~1 imply that for some~$a>0$, the couple
$(\gamma_n, a\|x\|^2)$
satisfies~($\tau$) for every~$n$.
However, the proof of Theorem~2 gives
a better constant~$a$ and is simpler.
\bigbreak

We will show now that Theorem~2 allows to recover the main
conclusion of the Gaussian concentration result of~\cfPi.

\proclaim Corollary 2. Let $\varphi$ be a 1-Lipschitz function on
$R^n$, and $X,Y$ two independent $n$-dimensional Gaussian vectors
with distribution equal to~$\gamma_n$. For every real number~$\lambda$
we have
$$ E e^{{\lambda\over{\sqrt 2}}(\varphi(X)-\varphi(Y))} \leq
e^{\lambda^2/2}.$$

\noindent{\bf Remark 3.} This inequality is optimal since there is
equality when $\varphi$ is a norm-one linear functional.
\smallbreak

\noindent Proof: Let $\psi = {{\lambda \varphi }\over \sqrt 2} \Carre w$,
where $\varphi$ is 1-Lipschitz on~$R^n$,
$w(y)={1\over 4}\|y\|^2$ and
$\lambda>0$. It is enough to apply~($\tau$) and notice
that $\psi(x) \ge
{{\lambda \varphi }\over \sqrt 2}(x)- \lambda^2/2$. Let $y$ be such that
$$ \psi(x) = {{\lambda \varphi }\over \sqrt 2}(y) + {1\over 4}
\|x-y\|^2 .$$
Then
$$ \psi(x) \ge {{\lambda \varphi }\over \sqrt 2}(x) -
{\lambda\over\sqrt 2}
\|x-y\| + {1\over 4}\|x-y\|^2 \ge {{\lambda \varphi }\over \sqrt 2}(x)
+ \min\{{1\over 4}u^2-{\lambda\over\sqrt 2}u;
u \in R  \}
= {\lambda\varphi\over\sqrt 2}(x) - \lambda^2/2.$$
\bigbreak

   The first part of the next Corollary is known (it is a
Poincar\'e-type inequality due to Chen~\cfCh).
\proclaim Corollary 3. If $\varphi$ is a Lipschitz function on~$R^n$,
we have
$$ {1\over 2} \int (\varphi(x)-\varphi(y))^2 d\gamma_n(x)
d\gamma_n(y) \le \int \|\nabla\varphi\|^2 d\gamma_n.$$
  More generally, this result holds for every
probability measure  $\mu$ on~$R^n$
such that $(\mu, w)$ satisfies~($\tau$) for a function
$w$ convex and greater than
${1\over 4}\|x\|^2$ in a neighborhood of $0$.

\noindent Proof: Let $u$ be a convex function such that
$u \le w$  and
$u(x)={1\over 4} \|x\|^2$ in a neighborhood of~$0$;
assume that
$\varphi$ is a compactly supported $C^1$-function. For~$t>0$
consider $\varphi_t=t\varphi$ and $\psi_t=\varphi_t\Carre u$.
One can check that
$$ \lim_{t\rightarrow 0} {{\psi_t(x)-\varphi_t(x)} \over t^2}
= -\|\nabla\varphi(x)\|^2$$
and the result follows easily from the property~($\tau$)
of $(\mu, u)$ applied to~$\varphi_t$, when
$t\rightarrow 0$.
\medbreak

\noindent{\bf Remark 4.} If $(\mu, w)$ satisfies~($\tau$) on~$R$,
where $\mu$ is such that
$ {1\over 2}\int (x-y)^2 d\mu(x) d\mu(y) > 1$,
we can apply Corollary~3 to
$\varphi(x)=x$ to conclude that $\{x; w(x) \ge
{1\over 4} x^2 \}$ is not a neighborhood of~$0$.
One can also show that if $(\mu, w)$ satisfies~($\tau$) on~$R$,
with $\mu$ symmetric
and $ \int x^2 d\mu(x) = 1 $, then $(\gamma,{1\over2}w"(0)t^2)$
also satisfies~($\tau$). This shows the necessity of a subquadratic
behavior at~$0$ for the function $w$.
\smallbreak

\proclaim Corollary 4. Let $\lambda_n$ denote the uniform probability
measure on $[0,1]^n$. There exists $a>0$ such that $(\lambda_n,
a\|x\|^2)$ satisfies~($\tau$) for every integer $n$
(one can take $a=\pi/2$).

\noindent Proof: Using Lemma~3, this follows from Theorem~2,
exactly like in Pisier~\cfPi.
\bigbreak

\noindent{\bf IV. Convex property ($\tau$)}
\medskip
 Assume that $w$ is a convex function on
some topological vector space~$X$
and that $\mu$ is a probability measure on~$X$. We say that the
couple $(\mu, w)$ satisfies the {\it convex} property~($\tau$) provided
$$ (\int e^{ \varphi\carre w} d\mu)(\int e^{ -\varphi} d\mu)
\leq 1$$
for every {\it convex} measurable function~$\varphi$ on~$X$.
\smallbreak

\proclaim Lemma 5. If $(\mu_i, w_i)$ satisfies
the convex property~($\tau$) on~$X_i$ for~$i=1,2$,
then $(\mu_1\otimes \mu_2, w)$ satisfies the convex property~($\tau$)
on~$X_1 \times X_2$, with
$$ w(x_1,x_2) = w_1(x_1) + w_2(x_2).$$

\noindent Proof: As in the proof of Lemma~1 we consider
$\varphi^y(x)=\varphi(x,y)$ for a convex function~$\varphi$
on $X_1 \times X_2$; we can apply the convex property~($\tau$) to
$\psi(y)=\log(\int e^{\varphi^y\carre w_1} d\mu_1)$
if we observe that $\psi$ is a convex function.
\medbreak

  We shall say that $\mu$ has diameter~$\le 1$ as a short way
to express that $\mu$ is supported by a set of diameter~$\le 1$.
The following Theorem is the equivalent in our language of
a result of Talagrand~\cfTab\ and its generalization by
Johnson and Schechtman~\cfJS.

\proclaim Theorem 3. Let $(X_i)$ be a family of normed spaces;
for each~$i$, let $\mu_i$ be a probability measure with
diameter~$\le 1$ on~$X_i$, and
$w_i(x) = {1 \over 4} \|x\|^2$ for $x \in X_i$. If $\mu$ is the
product of the family~$(\mu_i)$,
then $(\mu, w)$ satisfies the convex
property~($\tau$), with $w(x) = \sum_i w_i(x_i)$.

\noindent Proof: According to Lemma~5, we only need to prove
the result for a single probability measure~$\mu$ with
diameter~$\le 1$ on a normed space~$X$.
Let $A$ be a set of diameter~$\le 1$ that
supports~$\mu$, and let $\varphi$ be a convex function on~$X$;
assume without loss of generality
that $\inf\varphi(A) = 0$. Define $w(x) = {1 \over 4} \|x\|^2$
and $\psi = \varphi \Carre w$.
  Let $x\in A$, $\varepsilon > 0$ and $a\in A$
such that $\varphi(a) \le \varepsilon$.
We have, if $y = (1-\theta)x + \theta a$ and $0 \le \theta \le 1$
$$ \psi(x) \le \varphi(y) + {1\over 4} \|x-y\|^2
\le (1-\theta)\varphi(x) + \theta\varepsilon + {1\over 4}
\theta^2.$$
  Choosing an optimal~$\theta$,
we deduce from the above that
$\psi(x) \le k(\varphi(x))$ where $k(u)$ is equal to
$u-u^2$ if $0 \le u \le {1\over 2}$,
and to~${1 \over 4}$ if $u \ge {1 \over 2}$. We claim now that
$e^{k(u)} \le 2 - e^{-u}$. It follows that
$$ \int e^{\psi} d\mu \le 2 - \int e^{-\varphi} d\mu
\le (\int e^{-\varphi} d\mu)^{-1}$$
and this finishes the proof (the preceding computation was
inspired by \cfJS).
\smallskip
\noindent  Proof of the claim: For $0 \le u \le {1\over 2}$, we write
$${1\over 2} (e^{u-u^2} + e^{-u}) = e^{-u^2/2} \cosh(u-u^2/2)
\le e^{-u^2/2} \cosh(u) \le 1.$$
\smallbreak

\noindent{\bf Remark 5.} In the case of the probability~$\beta$
on~$[0,1]$
that gives measure $1\over 2$ to~$\{0\}$ and
$\{1\}$, it is easy to improve the function~$w$
from ${1 \over 4} x^2$ to ${1\over 2} x^2$.
\medbreak

\proclaim Corollary 5 \cfTab, \cfJS . Let $A$ be a measurable subset
of~$[0,1]^n$
and~$B$ its convex hull. For every product probability measure~$\mu$
on~$[0,1]^n$ we have
$$ \int e^{ {1\over 4} d_B^2} \, d\mu \le (\mu(A))^{-1}$$
where $d_B$ denotes the Euclidean distance to the set~$B$.
\bigbreak

\noindent{\bf References.}

\noindent\cfBo\ C.~Borell,
The Brunn-Minkowski inequality in Gauss space,
Inventiones Math. 30 (1975), 205--216.

\noindent\cfCh\ L.~Chen,
An inequality for the multivariate normal distribution,
J. Multivariate Anal. 12 (1982) 306--315.

\noindent\cfEh\ A.~Ehrhard,
Sym\'etrisation dans l'espace de Gauss,
Math. Scand. 53 (1983) 281--301.

\noindent\cfFLM\ T.~Figiel, J.~Lindenstrauss, V.~Milman,
The dimension of almost spherical sections of convex bodies,
Acta Math. 139 (1977) 53--94.

\noindent\cfJS\ W.~Johnson, G.~Schechtman,
Remarks on Talagrand's deviation inequality for Rademacher functions,
Texas Functional Analysis Seminar 1988--89.

\noindent\cfLe\ L.~Leindler,
On a certain converse of H\"older's inequality,
Acta Sci. Math. 33 (1972) 217--223.

\noindent\cfMi\ V.~Milman,
A new proof of the theorem of A.~Dvoretzky on sections of convex bodies,
Func. Anal. Appl. 5 (1971) 28--37.

\noindent\cfMS\ V.~Milman, G.~Schechtman,
Asymptotic theory of finite dimensional normed spaces,
Springer Lecture Notes in Math. 1200 (1986).

\noindent\cfPi\ G.~Pisier,
Probabilistic methods in the geometry of Banach spaces,
CIME Varenna 1985, Springer Lecture Notes in Math. 1206, 167--241.

\noindent\cfPib\ G.~Pisier,
Volume of convex bodies and Banach spaces geometry,
Cambridge University Press.

\noindent\cfPr\ A.~Prekopa,
On logarithmically concave measures and functions,
Acta Sci. Math. 34 (1973) 335--343.

\noindent\cfTa\ M.~Talagrand,
to appear.

\noindent\cfTab\ M.~Talagrand,
An isoperimetric theorem on the cube and the Khintchine-Kahane
inequalities,
Proc. Amer. Math. Soc. 104 (1988) 905--909.

\medskip\bigskip
\hbox to \hsize{ \sl \hfil
Universit\'e Paris 7, U.F.R. de Math\'ematiques }

\bye